\documentclass[10pt]{article}

\usepackage{amsmath}
\usepackage{amsmath,amsfonts,amssymb,amsthm,graphics,amscd}
\usepackage[lmargin=3cm, rmargin=3cm,tmargin=2.5cm,bmargin=2.5cm]{geometry}

\usepackage{hyperref}
\hypersetup{
    colorlinks=true,
    linkcolor=blue,
    filecolor=magenta,      
    urlcolor=cyan,
    pdftitle={Overleaf Example},
    pdfpagemode=FullScreen,
    }

\usepackage{color}

\newcommand{\bdf}{\begin{definition}}
\newcommand{\edf}{\end{definition}}
\newcommand{\bprop}{\begin{proposition}}
\newcommand{\eprop}{\end{proposition}}
\newcommand{\bcor}{\begin{corollary}}
\newcommand{\ecor}{\end{corollary}}
\newcommand{\bet}{\begin{theorem}}
\newcommand{\eet}{\end{theorem}}
\newcommand{\blm}{\begin{lemma}}
\newcommand{\elm}{\end{lemma}}
\newcommand{\bp}{\begin{proof}}
\newcommand{\ep}{\end{proof}}
\newcommand{\bex}{\begin{example}\rm}
\newcommand{\eex}{\end{example}}
\newcommand{\bexs}{\begin{examples}\rm}
\newcommand{\eexs}{\end{examples}}
\newcommand{\bremark}{\begin{remark}\rm}
\newcommand{\eremark}{\end{remark}}

\newtheorem{theorem}{Theorem}[section]
\newtheorem{lemma}[theorem]{Lemma}
\newtheorem{corollary}[theorem]{Corollary}

\newtheorem{example}[theorem]{Example}
\newtheorem{definition}[theorem]{Definition}
\newtheorem{proposition}[theorem]{Proposition}
\newtheorem{remark}[theorem]{Remark}

\newtheorem{examples}[theorem]{Examples}

\numberwithin{equation}{section}

\begin{document}

\title { On Left and Right semi-B-Fredholm operators}

\author { Alaa Hamdan and  
Mohammed Berkani\ }

\date{}

\maketitle \setcounter{page}{1}

\begin{abstract} To complete the study of Fredholm type operators of  \cite{HB1} and \cite{HB2},   we define in this paper the classes  of  left and right  semi-B-Fredholm operators (Definition \ref{left-right}). Then,  we prove  that an operator  $T \in L(X), X$ being a Banach space, is  a left( resp. right) semi-B-Fredholm if and only if $T$ is the direct sum of left(resp. right) semi-Fredholm operator and a nilpotent one. This result extend the earlier  characterization of  B-Fredholm operators as the direct sum of a Fredholm operator and a nilpotent one obtained in \cite[Theorem 2.7]{P7}.
 \end{abstract}

\renewcommand{\thefootnote}{}
\footnotetext{\hspace{-7pt} {\em 2020 Mathematics Subject
Classification\/}: 47A53,  16U90, 16U40.
\baselineskip=18pt\newline\indent {\em Keywords:\/}\,\,Left\,semi-B-Fredholm,\,right\,semi-B-Fredholm,\,Finite\,\,rank,\,Quasi-Fredholm,\,Drazin\,invertible,\,$p-$invertible. }

\section{Introduction}

Let $X$ be a Banach space,
let $L(X)$  be the Banach algebra of bounded linear operators acting on
the Banach space $X$ and let  $ T \in L(X).$
We will denote by $N(T)$  the null space of  $T$,  by $ \alpha(T)$
the nullity of $T$,  by $ R(T)$ the range of $T$ and by $\beta(T) $
its defect. If the range $ R(T)$ of $T$ is closed  and  $ \alpha(T) < \infty $
 (resp.  $\beta(T)< \infty $ ),
\noindent  then $T$ is called an upper semi-Fredholm (resp. a
lower semi-Fredholm) operator. A semi-Fredholm operator is an
upper or a lower semi-Fredholm operator. If both of $ \alpha(T) $
and $\beta(T) $ are finite, then $T $ is called a Fredholm operator
and the index of $T$  is defined by   $ ind(T) = \alpha(T) -
\beta(T). $  The notations $\Phi_+(X),  \Phi_-(X)$ and  $\Phi(X)$ will designate  respectively
the set of upper semi-Fredholm, lower semi-Fredholm and Fredholm operators. Define also  the sets:

\noindent 

$$\Phi_l(X) = \{T \in \Phi_+ (X) \mid \text { there exists a bounded projection of}\,  X\,  \text{onto}\, R(T)\} $$  and

$$\Phi_r (X) = \{T \in \Phi_-(X)\mid  \text { there exists a bounded projection of}\,  X\,  \text{onto}\, N(T) \}$$.

\noindent Recall that the Calkin algebra over $X$ is the Banach algebra, given by  the quotient algebra
$\mathcal{C}(X)=L(X)/K(X)$, where  $K(X)$ is the closed ideal of
compact operators on $X$.  Let $G_r$ and $G_l$   be the right and left, respectively, invertible elements
  of $ \mathcal{C}(X)$. From \cite[Theorem 4.3.2]{CPY}    and  \cite[Theorem 4.3.3]{CPY},
   it follows that $ \Phi_l (X)= \Pi^{-1}(G_l)$  and  $ \Phi_r (X)= \Pi^{-1}(G_r),$
     where $ \Pi: L(X) \rightarrow  \mathcal{C}(X)$ is the natural projection. We observe that $\Phi(X)= \Phi_l (X)\cap \Phi_r(X).$

  \bdf \label{left-right semi}The elements of  $ \Phi_l (X)$ and   $\Phi_r (X)$ will be  called respectively left
    semi-Fredholm operators and   right  semi-Fredholm operators .

  \edf

In 1958, in his paper \cite{DR},  the author  extended the
concept of invertibility in  rings and semigroups and introduced a
new kind of inverse, known now as the  Drazin inverse.

\bdf \label{def-Drazin}An element  $a$  of  a semigroup $ \mathcal{S}$ is called
Drazin invertible if there exists an element $b \in \mathcal{S}$ written $ b= a^d$  and called the Drazin inverse of $a,$
satisfying the following equations

\begin{equation}\label{Drazin}
  ab = ba, b=ab^2, a^k= a^{k+1}b,
 \end{equation}

 for some nonnegative integer $k.$ The least nonnegative
integer $k$ for which these equations holds is called  the Drazin index $
i(a)$ of $a.$

\edf

\noindent It follows from \cite{DR} that a Drazin invertible
element in a semigroup has a unique Drazin inverse.\\
If $A$ is a unital Banach algebra,  it is easily seen that the conditions  (\ref{Drazin}) for the Drazin invertibility of an element $ a\in A$  are equivalent to the existence of an idempotent $ p \in A $ such that  $$ap=pa, ap \text{ is nilpotent and} \,\,  a+p   \text{ is   invertible}. $$

\noindent In 1996, in  \cite[Definition 2.3]{Koliha}, the author extended
the notion of  Drazin invertibility as follows.

\bdf An element  $a$  of a Banach algebra $A$ will be said to be
generalized Drazin invertible if there exists $b \in A $ such that
$ bab=b, ab=ba$ and  $ aba-a$ is a quasinilpotent element in $A$.
\edf

In \cite[Theorem 4.2]{Koliha}, the author proved that  $a\in A$ is generalized Drazin
invertible if and only if  there exists $ \epsilon >0,$ such that
for all $\lambda \in \mathbb{C}$ such that $ 0 < \mid \lambda \mid <  \epsilon,
$ the element  $a-\lambda e $ is invertible and he proved
 that a generalized Drazin
invertible element  has a unique generalized Drazin
inverse. He also proved that an element  $a\in A$ is generalized Drazin
invertible if and only if  there exists an idempotent $p \in A$  commuting with $a,$
such that $ a+p$ is invertible in $A$ and $ap$ is quasinilpotent. \\
Recall that in a ring $A$ with a unit, for $a \in A,$ the commutant  $comm(a)= \{ x \in A \mid xa= ax \}$   and the bicommutant $ comm^2(a)=\{ x \in A \mid xy= yx$ for all  $y \in comm(a)\}.$ 
 
 Let $A$ be a ring with a unit and let $p$ be an idempotent element in $A$.  In  \cite[Definition 2.2]{P50},
the concepts of left $p-$invertibility, right $p-$invertibility and $p-$invertibility where defined  as follows.

\bdf Let $ a \in A.$   We will say that

\begin{enumerate}

\item $a$ is left  $p$-invertible if $ ap=pa$ and  $ a+p$ is left
invertible in $comm(p).$

\item  $a$ is right  $p$-invertible if $ ap=pa$ and $ a+p$ is right
invertible in $comm(p).$

 \item $a$ is
$p$-invertible if $ ap=pa$    and $ a+p$ is invertible.

\end{enumerate}

\edf

Moreover in \cite[Definition 2.11]{P50}, left and right Drazin invertibility where defined as follows.

\bdf \label{Drazin $p$-invertibility}  We will say that an element
$a  \in A $ is left Drazin invertible (respectively right Drazin
invertible) if there exists an idempotent $p \in A$ such that $a$
is left $p$-invertible (respectively right $p$-invertible) and
$ap$ is nilpotent.

\edf

For  $T \in L(X),$  we will say that  a subspace $M$
 of $X$  is {\em $T$-invariant} if $T(M) \subset M.$ We define
$T_{\mid M}:M \to M$ as $T_{\mid M}(x)=T(x), \, x \in M$.  If $M$
and $N$ are two closed $T$-invariant subspaces of $X$ such that
$X=M \oplus N$, we say that $T$ is {\em completely reduced} by the
pair $(M,N)$ and it is denoted by $(M,N) \in Red(T)$. In this case
we write $T=T_{\mid M} \oplus T_{\mid N}$ and we say that $T$ is a
{\em direct sum} of $T_{\mid M}$ and $T_{\mid N}$.

It is said that $T \in L(X)$ admits a  generalized
Kato decomposition, abbreviated as GKD,  if there exists $(M, N)
\in Red(T)$ such that $ T_{\mid M}$ is Kato and $T_{\mid N}$ is
quasinilpotent. Recall that an operator $T \in L(X)$ is a {\em Kato operator}
if $R(T)$ is closed and $Ker(T) \subset R(T^n)$ for every $ n \in
\mathbb{N}$.
\bremark For $T \in L(X),$ to say that a pair $(M,N)$ of closed subspaces of $X$ is  in $Red(T),$
means simply that there exists an idempotent $P \in L(X)$ commuting with $T.$ Indeed if $(M,N) \in Red(T),$
let $P$ be the projection on $M$ in parallel to $N.$ Then since $(M,N) \in Red(T),$ we see that  $P$ commutes  with $T.$
Conversely, given an idempotent $P \in L(X)$ commuting with $T,$ if we set $M= P(X)$ and $N= (I-P)X,$ then it is clear that $(M,N) \in Red(T).$

\eremark

For $T \in L(X)$ and a nonnegative integer $n,$ define $T_{[n]}$ to
be the restriction of $T$ to $R(T^n)$ viewed as a map from $R(T^n)$
into $R(T^n)$ (in particular $T_{[0]}=T).$ If for some integer $n$
the range space $R(T^n)$ is closed and $T_{[n]}$ is an upper (resp.
a lower) semi-Fredholm operator, then $T$ is called an  upper (resp.
a lower) semi-B-Fredholm operator.
Moreover, if $T_{[n]}$ is a Fredholm operator, then $T$ is called a
B-Fredholm operator (see \cite{P7} for more details). The following theorem gives a characterization  of B-Fredholm operators 

\bet \cite[Theorem 2.7]{P7} \label{decomposition} Let   $ T  \in L(X). $ Then  T  is   a
 B-Fredholm  operator if and only if  there exists two closed subspaces  $(M, N)$  in $Red(T)$  such that:
 
 \begin{enumerate}
\item $T(N)\subset N$ and  $T_{\mid N}$ is a  nilpotent operator,
\item $T(M) \subset M$  and  $T_{\mid M }$ is a Fredholm  operator. 
\end{enumerate}
\eet


\bdf  \cite[Definition 2]{HA}  An element $a$  of a   ring A with
a unit $e$ is quasinilpotent if, for every x commuting with a, $e -xa $ is invertible in $A.$ \edf


\bdf \label{gdrazin}  \cite[Definition 3.3]{HB2} Let  $T \in L(X).$
\begin{enumerate}

\item  We will say that $\pi(T)$ is generalized Drazin invertible (resp. left generalized Drazin invertible, resp. right  generalized Drazin invertible) in $\mathcal{C}_0(X),$ if there exists an idempotent $P \in L(X) $   such that  $\pi(P)\in comm^2(\pi(T)),$   $ \pi(P)\pi(T)$ is quasi-nilpotent  and $\pi(T)+ \pi(P)$ is  invertible (resp. $\pi(T)$  left $\pi(P)-$invertible, resp. $\pi(T)$ right $\pi(P)-$ invertible) in $\mathcal{C}_0(X).$

\item We will say that $\Pi(T)$ is generalized Drazin invertible (resp. left generalized Drazin invertible, resp. right  generalized Drazin invertible) in $\mathcal{C}(X),$ if there exists an idempotent $P \in L(X) $   such that  $\Pi(P)\in comm(\Pi(T)),$   $ \Pi(P)\Pi(T)$ is quasi-nilpotent  and $\Pi(T)+ \Pi(T) $ is  invertible (resp. $\Pi(T)$  left $\Pi(P)-$invertible, resp. $\Pi(T)$ right $\Pi(P)-$ invertible) in $\mathcal{C}(X).$ 

\end{enumerate}

\edf

As mentionned before,  generalized Drazin invertibility in Banach algebras was introduced in \cite{Koliha}, while  generalized Drazin invertibility in rings  was introduced in \cite{KP}.

 First let us recall the following important result, which links the idempotents of the Calkin algebra   $ \mathcal{C}(X)$  to idempotents of $L(X).$ 
 
 \blm \cite[Lemma 1]{BAR} \label{Barnes}
Let $p$ be an idempotent element of the Calkin algebra $ \mathcal{C}(X).$  Then there exists
an idempotent $P \in L(X)$ such that $ \Pi(P)= p.$
\elm
 
  In \cite[Theorem 2.1]{HB2}, we extended Lemma \ref{Barnes} to the case of the algebra $ \mathcal{C}_0(X),$  where for a seek of notations simplification we set 
 $ \mathcal{C}_0(X)= L(X)/F_0(X),  F_0(X)$ being the ideal of finite rank operators. Then $\pi: L(X) \rightarrow \mathcal{C}_0(X)$ is the natural projection.

 \blm \label{idempotent} \cite[Theorem 2.1]{HB2}
Let $p_0$ be an idempotent element of the algebra $ \mathcal{C}_0(X).$  Then there exists
an idempotent $P \in L(X)$ such that $ \pi(P)= p_0.$
\elm
 As this paper is a continuation of \cite{HB1} and \cite{HB2}, we begin by giving a brief  summary of the results of \cite{HB1} and \cite{HB2}. So letting $p=\Pi(P)$ being  an idempotent  in the Calkin algebra and    using the results on $p$-invertibility obtained  in  \cite{P50},   we introduced  the class $\Phi_P(X)$ of $P$-Fredholm (respectively the class $\Phi_{lP}(X)$ of left semi-$P$-Fredholm and the class $\Phi_{rP}(X)$ of right  semi-$P$-Fredholm)  operators, in a similar way as  the corresponding classes of Fredholm, left semi-Fredholm and right semi-Fredholm operators. Then we observed that
$\Phi_P(X)= \Phi_{lP}(X)\cap \Phi_{rP}(X).$

 Moreover, using left and right generalized Drazin invertibility in the Calkin algebra, we  introduced the classes $\Phi_{l\mathcal{P}B}(X) $ and $\Phi_{r\mathcal{P}B}(X) $ of left and right  pseudo semi-B-Fredholm operators, completing in this way the study of the class  $\Phi_{\mathcal{PB}}(X) $ of pseudo B-Fredholm operators inaugurated  in \cite{P45} and  proving that   $\Phi_{\mathcal{P}B}(X)= \Phi_{l\mathcal{P}B}(X)\cap \Phi_{r\mathcal{P}B}(X).$

 Based on left and right  Drazin invertibility in the Calkin algebra, we
 introduced  the classes $\Phi_{l\mathcal{W}B}(X),$  $\Phi_{r\mathcal{W}B}(X) $  of  left and right  weak semi-B-Fredholm operators, completing in this way the study of the class  $\Phi_{\mathcal{W}B}(X) $ of weak B-Fredholm operators inaugurated  in \cite{P45} and proving that  $\Phi_{\mathcal{W}B}(X)= \Phi_{l\mathcal{W}B}(X)\cap \Phi_{r\mathcal{W}B}(X).$

Though, the weak B-Fredholm operators and  pseudo B-Fredholm  operators, where not  explicitly defined in \cite{BO}, they where characterized there  by  Drazin (generalized Drazin) invertible elements in the Calkin algebra in \cite[Theorem 6.1,ii]{BO}  and  \cite[Theorem 6.1, i]{BO} respectively.

 In summary, the study in \cite{HB1} was based on different type of invertibility in the Calkin algebra.  So  the classes of  Fredholm type operators of \cite{HB1} where defined   \textit{modulo} the Calkin algebra.  As seen in \cite{HB1},  these classes of operators obeys   the following strict inclusions relations, $$ \Phi(X) \subsetneq \Phi_B(X) \subsetneq \Phi_{\mathcal{W}B}(X) \subsetneq \Phi_{\mathcal{P}B}(X)\subsetneq \Phi_\mathbb{P}(X)$$  where  $ \Phi_{\mathbb{P}}(X)=\bigcup\limits_{ \{ P\in L(X) \mid P^2=P\}}\Phi_P(X). $

In a similar way, in \cite{HB2},  the classes of P-Fredholm operators and pseudo-B-Fredholm operators where defined  \textit{modulo} the algebra  $\mathcal{C}_0(X).$ 
 
The aim of the present paper is to complete the study of Fredholm type operators done in  \cite{HB1} and \cite{HB2}.  Here  we define  the classes  of  left and right  semi-B-Fredholm  operators (Definition \ref{left-right}) and we prove  that $T$ is  a left(resp. right) semi-B-Fredholm operator if and only if $ T$ is the direct sum of left (resp. right) semi-Fredholm operator and a nilpotent one. This result extend the earlier  characterization of  B-Fredholm operators as the direct sum of a Fredholm operator and a nilpotent one obtained in \cite[Theorem 2.7]{P7}. As a consequence, an operator $T$ is a B-Fredholm operator  if and only if $ \pi(T)$ is left and right Drazin invertible in the algebra $\mathcal{C}_0(X),$ where  $\pi: L(X) \rightarrow \mathcal{C}_0(X)$ is the natural projection. As a corollary, it follows immediately that left and right semi-B-Fredholm operators are stable under finite rank perturbation.

\bremark \label{not1} Unless mentioned otherwise, in all this paper,  we will use the following notations.   For  $T \in L(X)$ and  $ P \in L(X)$ an idempotent element of $L(X),$ we will set  $X_1= R(P), X_2= N(P)= R(I-P),$   $T_1= (PTP)_{\mid_{X_1}}$ and  $T_2= (I-P)T(I-P)_{\mid_{X_2}}.$

  $\bullet$  If $P$ commutes with $T,$ then  we have  $T= T_1 \oplus T_2,$  here $T_1= T_{\mid_{X_1}}$ and  $T_2= T_{\mid_{X_2}}.$ In this case $(X_1, X_2) \in Red(T).$

   $\bullet$ If $\pi(P)$ commutes with $\pi(T)$ in  $\mathcal{C}_0(X),$   then we have $T= TP + T(I-P)= PTP + (I-P)T(I-P) + F,$ where $F \in L(X)$ is a finite rank operator.
   So  $T= T_1 \oplus T_2 + F.$ In this case $(X_1, X_2) \in Red(PTP)$ and $(X_1, X_2) \in Red((I-P)T(I-P)).$

   $\bullet$ $ I_1= P_{\mid X_1},   I_2= (I-P)_{\mid X_2}.  $

   $\bullet$  When we use the homomorphisms $\pi$,  we  mean the natural projection from $L(E)$ onto the algebra $\mathcal{C}_0(E),$  $E$  being a Banach subspace of $X$  appearing in the context of the use of $\pi.$

\eremark

\section{Quasi-Fredholm Operators}

\bdf \cite{LAB} Let   $ T \in L(X) $  and let
$$\Delta (T)=\{n\in \mathbb{N}: \;\;\forall m\in \mathbb{N}\,\ m \geq n \Rightarrow
(R(T^n)\cap N(T))\subset (R(T^m)\cap N(T))\}.$$
Then the  degree of  stable iteration $ dis(T) $ of   $T$  is defined as
$dis(T)=inf\Delta (T)$ (with  \\
\noindent dis(T) = $ \infty $  if $\Delta(T) = \emptyset$). \\
\edf 
\bdf \cite{LAB} \label{QF1} Let  $ T \in L(X). $   Then  $T$  is called a quasi-Fredholm
 operator of degree d  if  there is an integer  $ d\in \mathbb{N}$ such that: \\
\noindent (a) $dis(T)=d,$\\
(b) $ N(T)\cap R(T^d)$  is  a closed and complemented  subspace of  $X.$\\
(c) $R(T)+N(T^d)$  is a closed and complemented subspace of  $X.$ \\
In the sequel, the symbol  $ QF(d) $ will denote the set of quasi-Fredholm operators  of degree d.
\edf

Recall that an operator $T \in L(X)$ is called a regular operator if $ N(T)$ and $R(T)$ are complemented subspaces of $X.$ See \cite[C.12.1]{MU} and \cite[Propositon 1]{MU} for more details. 

\bet \cite[Th\'eor\'eme 3.2.2]{LAB} \label{LAB1}
Let $T \in L(X).$ Then $T$ is a   quasi-Fredholm operator of degree $d, d \geq 1,$  if and only if there exists two closed
subspaces  $M,  N$  of $ X $ and an integer  $ d \in \bf{N} $ such that   $ X = M\oplus N $  and:

\begin{enumerate}
\item $T(N)\subset N $  and  $ T_{\mid N}$ is a  nilpotent operator of degree $d$,
\item $ T(M) \subset M,  
 N(T_{\mid M }) \subset\underset{m}\cap R( (T_{\mid M })^m)  $ and $T_{\mid M }$ is a regular operator.

\end{enumerate}
\eet

\bremark 

\begin{enumerate}
\item Tough Definition \ref{QF1} and Theorem \ref{LAB1} had been considered in \cite{LAB} in the case of Hilbert spaces, it is mentioned in \cite[Remarque, page 206]{LAB}
that it holds also in the case of Banach spaces, with the same proof. In fact the same proof exactly  holds in the case of Banach spaces, that's why  it was not included in \cite[Theorem 2.7]{P7} and nor included  here. See also the comments in \cite [C.22.5]{MU}.

\item In \cite{ZD}, the authors gave an extensive study of the class of quasi-Fredholm operators under the appellation of ``Saphar type Operators".

\end{enumerate}

\eremark

\bdf \cite{MU} \label{QF2} Let  $ T \in L(X). $   Then  $T$  is called a weak quasi-Fredholm
 operator of degree d  if  there is an integer   $d\in \mathbb{N}$ such that $dis(T)=d $  and  $R(T^{d+1})$   is closed.\\

In the sequel, the symbol  $ wQF(d) $ will denote the set of weak quasi-Fredholm operators  of degree d.
\edf

\bremark \label{wQF(d)}

\begin{enumerate}

\item Every quasi-Fredholm operator is a weak quasi-Fredholm operator, but the converse may not hold. If $T$ is quasi-Freddholm operator of degree $d,$ then following the same proof as in \cite[Proposition 3.3.2]{LAB} (Proven in the case of Hilbert spaces, but the same proof holds also in Banach spaces), we can prove that  $ R(T^{d+1})$ is closed.
    
     In the case of Hilbert spaces, the two classes coincides, that's  $ QF(d)= wQF(d).$

\item Weak quasi-Fredholm operators  are called in \cite{MU} and \cite{KMMP}  ``Quasi-fredholm operators".  We  call them here weak quasi-Fredholm operators  to avoid confusion with the class of quasi-Fredholm operators of Definition \ref{QF1}. For the purpose of our paper, the definition given in  Definition \ref{QF2} is in fact an equivalent definition to that  used in \cite{MU} and \cite{KMMP}.
\item It follows from \cite[Lemma 17]{MU} that if $T$ is weak quasi-Fredholm with $dis(T)= d,$ then $R(T^n)$ is closed for every $ n\geq d.$

\item If $\in QF(d),$ then  from \cite[Lemma 17]{MU}, the subspaces  $N(T)\cap R(T^d)$   and  $R(T)+N(T^d)$ are closed subspaces of  $X.$\\

\end{enumerate}
\eremark

We will use the following notations from \cite[p. 181]{MU}.  
For closed  subspaces M, L of Banach space X, we write $ M \overset{e}{\subset} L $  ($M$ is essentially contained in $L$) if there exists a finite-dimensional subspace $ F \subset X,$ such that $ M \subset L+F.$ Equivalently $dim( M/(M\cap L)= dim ( M+L/M) < \infty. $  We will write $M \overset{e}{=} L $ if  $M \overset{e}{\subset} L$ and $L \overset{e}{\subset} M.$

\bet \label{perturbation} Let $T \in L(X)$  be a quasi-Fredholm operators and let $ F\in L(X)$ be a finite rank operator, then $T+F$ is a quasi-Fredholm operator. Moreover, if $d= dis(T)$ and
 $d'= dis(T+F),$ the following properties  holds.

\begin{enumerate}
\item  If $ N(T) \cap R(T^d)$  is of finite dimension, then  $N(T+F) \cap R((T+F)^{d'}) $ is also of finite dimension.

\item  If $R(T)+N(T^d)$  is of finite codimension, then  $R(T+F) + N((T+F)^{d'})$ is also of finite codimension.

\end{enumerate}

\eet

\bp As every quasi-Fredholm operator is a weak quasi-Fredholm, then from \cite[C.22.4, Table 2]{MU}, it follows that $T+F$ is also a weak quasi-Fredholm operator. 
For $n$ large enough, using \cite[Theorem]{KMMP}  we have $N(T+F) \cap R((T+F)^{d'})= N(T+F) \cap R((T+F)^{n}) \overset{e}{=}   N(T) \cap R(T^n)= N(T) \cap R(T^d).$  As  $N(T) \cap R(T^d)$ is a  complemented subspace of $X,$  then from \cite[Appendix 1, Theorem 25, iii)]{MU}, it follows that $N(T+F) \cap R((T+F)^{d'})$ is also complemented.

Now as  $ T \in wQF(d),$ then from \cite[C.22.6, p. 217]{MU}, its adjoint $T^* \in wQF(d)$   in $X^*.$ As $F$ is finite dimensional, then $F^*$ is also finite dimensional and from \cite[C.22.4, Table 2]{MU}, $T^*+F^*$ is also a weak quasi-Fredholm operator. We have  $ d'= dis(T^* +F^*)$  and for  $n > d +d',$  as $R(T)+N(T^d)$ is  complemented, then  $R(T)+N(T^n)$  is also complemented. So there  exists a closed subspace $E$ of $X$ such that $ X=[R(T)+N(T^n)] \oplus E.$ From \cite[Lemma 17]{MU}, $R(T^{n})$ is closed   and so $$X^*= (R(T)+N(T^n))^{\bot} \oplus E^{\bot} =[N(T^*) \cap R((T^*)^{n})] \oplus E^{\bot}= [N(T^*) \cap R((T^*)^{d'})] \oplus E^{\bot}.$$  From \cite{KMMP}, we have  $ N(T^*+F^*) \cap R((T^*+F^*)^{d'}) \overset{e}{=} N(T^*) \cap R((T^*)^{d}).$ So there exists two finite dimensional subspaces $ H_1, H_2$ of $X^*$ such that $$   N(T^*) \cap R((T^*)^{d}) \subset  N(T^*+F^*) \cap R((T^*+F^*)^{d'})  + H_1 $$ and  $$    N(T^*+F^*) \cap R((T^*+F^*)^{d'})  \subset  N(T^*) \cap R((T^*)^{d})   + H_2. $$

Let $ P: X^* \rightarrow X^*$ be the projection on $ N(T^*) \cap R((T^*)^{d})$ in parallel to $E^{\bot}.$  Then $$ N(T^*+F^*) \cap R((T^*+F^*)^{d'})  \subset  N(T^*) \cap R((T^*)^{d})   +(I-P) H_2. $$ 

and

 $$X^*=  N(T^*) \cap R((T^*)^{d})  \oplus E^{\bot}  \subset N(T^*+F^*) \cap R((T^*+F^*)^{d'})   + H_1 + E^{\bot}\subset X^*. $$

Thus $ N(T^*+F^*)\cap R((T^*+F^*)^{d'})  + E^{\bot}$  is of finite codimension in $X^*.$   So $(R(T+F) + N((T+F)^{d'}))^{\bot}+ E^{\bot} $ is closed.  
As $R(T+F) + N((T+F)^{d'})$ is closed (because  $T+F$ is a weak quasi-Fredholm operator)   and $E$ is closed,  then from \cite[Theorem 13]{MU}  
$R(T+F) + N((T+F)^{d'})+E$ is closed. 

Let $K$ be  a finite dimensional subspace of $X^*$ such that  
\begin{equation}\label{equality1}
X^* =[ [ N(T^*+F^*) \cap R((T^*+F^*)^{d'})] + E^{\bot}] \oplus K. 
\end{equation}

 Set $ m=dim(K),$ then  there exists  $ \{g_1, g_2,..., g_m\} \subset X^*$  such that 
 $K= vect\{g_1, g_2,..., g_m\}= (\underset{i=1}{\overset{m}{\bigcap}} Ker(g_i))^{\bot},$  where $vect\{g_1, g_2,..., g_m\}$ is the vector subspace of $X^*$  generated by the set $\{g_1, g_2,..., g_m\}.$  So $$ X^* =[ [ N(T^*+F^*) \cap R((T^*+F^*)^{d'})] + E^{\bot}] \oplus (\underset{i=1}{\overset{m}{\bigcap}} Ker(g_i))^{\bot}.$$

We have $$[ N(T^*+F^*) \cap R((T^*+F^*)^{d'})] \cap E^{\bot} \subset [N(T^*) \cap R((T^*)^{d})   + (I-P)(H_2)]\cap E^{\bot} $$ and so  $$[ N(T^*+F^*) \cap R((T^*+F^*)^{d'})] \cap E^{\bot}  \subset  (I-P)(H_2). $$
Hence $[ N(T^*+F^*) \cap R((T^*+F^*)^{d'})] \cap E^{\bot}$ is of finite dimension. Let $n$ be its dimension,  then  there exists  $ \{f_1, f_2,..., f_n\} \subset E^{\bot}$  such that 

\begin{equation} \label{equality2}
[ N(T^*+F^*) \cap R((T^*+F^*)^{d'})] \cap E^{\bot}= vect\{f_1, f_2,..., f_n\}= (\underset{i=1}{\overset{n}{\cap}} Ker(f_i))^{\bot}. 
\end{equation}

  Let $ \{x_1, x_2, ..., x_n\} \subset X $ such that  $ (\underset{i=1}{\overset{n}{\cap}} Ker(f_i)) \oplus vect\{x_1, x_2, ..., x_n\} = X.$ 
Then 

$$ E^{\bot} = vect\{f_1, f_2,..., f_n\} \oplus [E \oplus vect\{x_1, x_2, ..., x_n\}]^{\bot}= (\underset{i=1}{\overset{n}{\cap}} Ker(f_i))^{\bot} \oplus G^{\bot} $$ 
  where  $ G= E \oplus vect\{x_1, x_2, ..., x_n\}.$ Now using \ref{equality2}, we obtain  $$ [ N(T^*+F^*) \cap R((T^*+F^*)^{d'})] \cap G^{\bot}=    [ N(T^*+F^*) \cap R((T^*+F^*)^{d'})] \cap E^{\bot} \cap G^{\bot}\\$$$$=  (\underset{i=1}{\overset{n}{\cap}} Ker(f_i))^{\bot} \cap G^{\bot} = \{0\}. $$   With the aid of   \ref{equality1}, we obtain
  
  $$ X^* =[ [ N(T^*+F^*) \cap R((T^*+F^*)^{d'})] + E^{\bot}] \oplus K= [ N(T^*+F^*) \cap R((T^*+F^*)^{d'})] \oplus G^{\bot} \oplus K.$$

Then  $$ X^* = [ N(T^*+F^*) \cap R((T^*+F^*)^{d'})] \oplus [(E + vect\{x_1, x_2, ..., x_n\}) \cap (\underset{i=1}{\overset{m}{\cap}} Ker(g_i)\})]^{\bot}$$ 
 $$= [R(T+F) +N((T+F)^{d'})]^{\bot}  \oplus  [(E + vect\{x_1, x_2, ..., x_n) )\cap [(\underset{i=1}{\overset{m}{\cap}} Ker(g_i)]^{\bot} $$

As $R(T+F) +N((T+F)^{d'})$ and $(E + vect\{x_1, x_2, ..., x_n) )\cap (\underset{i=1}{\overset{m}{\cap}} Ker(g_i)$  are closed and the sum of their annihilators is closed, then from \cite[Theorem 13]{MU}, we have 
 $$ X= [R(T+F) +N((T+F)^{d'})]  \oplus  [(E + vect\{x_1, x_2, ..., x_n\}) \cap (\underset{i=1}{\overset{m}{\cap}} Ker(g_i)\})] $$
and
$[R(T+F) +N((T+F)^{d'})]$ is complemented in $X.$ Finally we see that  $T+F$ is a quasi-Fredholm operator.

 Now if $N(T)\cap R(T^d)$  is of finite dimension,  then $T$ is an upper semi-B-Fredholm operator. From \cite[Proposition 2.7]{P11},  $T+F$ is also an upper semi-B-Fredholm operator.  Then
 $ N(T+F)\cap R((T+F)^{d'}) $ is of finite dimension.  

Similarly if  $R(T)+N(T^d)$ is of finite codimension, then $ {\frac { R(T^{d})}{R(T^{d+1})}} $ is of finite dimension because it is isomorphic to   $\frac{X}{R(T) + N(T^{d})}.$  So $T$ is a lower semi-B-Fredholm operator. 
From \cite[Proposition 2.7]{P11},  $T+F$ is also a lower semi-B-Fredholm operator. Hence, ${\frac { R((T+F)^{d'})}{R((T+F)^{d'+1})}}$  is of finite dimension  and so 
$R(T+F) + N((T+F)^{d'})$ is also of finite codimension, because $\frac{X}{R(T+F) + N((T+F)^{d'})}$   is isomorphic to ${\frac { R((T+F)^{d'})}{R((T+F)^{d'+1})}}$    \ep 

\bexs \label{examples}
\begin{enumerate}

\item  Let $T \in L(X)$ be a left semi-Fredholm operator. Then $T$ is a quasi-Fredholm operator. Indeed as $N(T)$ is of finite dimension, then the sequence  $( N(T) \cap R(T^n))_n$ is a stationary sequence. Hence $ d=dis(T)$  is finite and $ N(T)\cap R(T^d)$ is of finite dimension and so it is a complemented subspace of $X.$ Moreover, as $T$ is left semi-Fredholm operator, then $T^d$ is also a left semi-Fredholm operator. So $ R(T^d)$ is complemented in $X$ and  it is  also the case of $ N(T) + R(T^d),$ because $N(T)$ is of finite dimension. Therefore $T$ is a quasi-Fredholm operator.
  
\item   Similarly if  $T \in L(X)$ is a right semi-Fredholm operator, then $R(T)$ is of finite codimension.   Then the sequence  $(R(T) + N(T^n))_n$ is an increasing sequence of subspaces of finite codimension. So it is a stationary sequence and from \cite[Lemma 3.5 ]{KAS},  we have: $$ \frac{N(T) \cap R(T^{n})} {N(T) \cap R(T^{n+1})}\simeq
         \frac{ N(T^{n+1}) + R(T)}  {N(T^{n}) + R(T)} $$
      
      Where the symbol $\simeq$ stands for isomorphic. Hence the sequence $(N(T) \cap R(T^{n}))_n$ is a stationary sequence,   $ d=dis(T)$  is finite, $R(T) + N(T^d)$ is of finite codimension and so it is a complemented subspace of $X.$
     Now, as $T$ is right semi-Fredholm, then $T^{d+1}$ is also a right semi-Fredholm. So there exists a closed subspace $E$ of $X$ such that $ X= N(T^{d+1}) \oplus E.$ Then we can easily verify  that $(N(T)\cap R(T^d)) \cap T^d(E)= \{0\}$ and $ R(T^d)= N(T)\cap R(T^d) \oplus T^d(E).$    As the sum  and the intersection of the subspaces  $N(T)\cap R(T^d)$ and  $T^d(E)$ are closed, then from Neubauer Lemma \cite[C.20.4]{MU}, $T^d(E)$ is closed. Thus $N(T)\cap R(T^d)$ is complemented in $R(T^d).$ As  $R(T^d)$ is of finite codimension in $X$, then $N(T)\cap R(T^d)$ is complemented in $X.$


\end{enumerate}

\eexs

 \section{Left and Right semi-B-Fredholm operators}

In addition to the decomposition theorem \ref{decomposition}, we know from \cite[Theorem 3.4]{P10} that  $T \in L(X)$ is a B-Fredholm operator
if and only if $\pi(T)$ is Drazin invertible in the algebra $\mathcal{C}_0(X).$  This result and Definition \ref{left-right semi} motivate the following definition of left and right semi-B-Fredholm operators.

\bdf \label{left-right} Let $ T \in L(X).$    We will say that:

\begin{enumerate}
  \item $T$ is a left  semi-B-Fredholm operator if $ \pi(T)  $ is left  Drazin  invertible in $\mathcal{C}_0(X).$
  \item $T$ is a right semi-B-Fredholm operator if $ \pi(T)$ is right Drazin invertible in $\mathcal{C}_0(X).$
  \item $T$ is a  strong  semi-B-Fredholm operator if it is left  or right   semi-B-Fredholm operator.

\end{enumerate}

\edf

The study of  strong semi- B-Fredholm operators involves the following classes of operators.

\bdf \label{power finite rank}

\begin{enumerate}

  \item $T \in L(X)$  is called a power finite rank\textbf{\Large -}left semi-Fredholm operator if
 there exists $(M, N) \in Red(T)$ such that $
T_{\mid M}$ is a power finite rank   operator and $T_{\mid N} \in \Phi_l(N) .$

  \item $T \in L(X)$  is called  a power finite rank\textbf{\Large -}right semi-Fredholm operator if
 there exists $(M, N) \in Red(T)$ such that $ T_{\mid M}$ is a power finite rank-right  operator and $T_{\mid N} \in \Phi_r(N) .$

\end{enumerate}

\edf

To characterize left and right semi-B-Fredholm operators, we will use the following classes of operators, which are subclasses of the corresponding  classes of Definition \ref{power finite rank}.

\bdf

\begin{enumerate}
  \item $T \in L(X)$  is called a nilpotent\textbf{\Large-}left semi-Fredholm operator if
 there exists $(M, N) \in Red(T)$ such that $
T_{\mid M}$ is a nilpotent   operator and $T_{\mid N} \in \Phi_l(N) .$
  \item  $T \in L(X)$  is called a nilpotent\textbf{\Large-}right  semi-Fredholm operator if
 there exists $(M, N) \in Red(T)$ such that $
T_{\mid M}$ is nilpotent  operator and $T_{\mid N} \in \Phi_r(N) .$

  \item  $T \in L(X)$  is called a nilpotent\textbf{\Large-}Fredholm operator if
 there exists $(M, N) \in Red(T)$ such that $
T_{\mid M}$ is nilpotent  operator and $T_{\mid N}$ is a Fredholm operator.
\end{enumerate}

\edf

\bremark 
In \cite{ZD}, the authors gave an extensive study of the classes of ``essentially left (resp. right) Drazin invertible operators" corresponding to the classes   of nilpotent\textbf{\Large-}left (resp. right)  semi-Fredholm operators we are considering here.

\eremark

In the next theorem, we prove the equivalence of Drazin invertibility with left and right Drazin invertibility for elements of $ \mathcal{C}_0(X).$

\bet \label{first equivalence} Let $T \in L(X),$   then the following properties
are equivalent.

\begin{enumerate}

\item $\pi(T)$ is  Drazin invertible in $\mathcal{C}_0(X).$

\item $\pi(T)$  is left and  right Drazin invertible in $\mathcal{C}_0(X).$

\end{enumerate}
\eet

\bp Since the implication $ 1) \Rightarrow 2)$ is trivial, we only need to prove the implication $ 2) \Rightarrow 1).$

So assume that $\pi(T)$  is left and  right  Drazin invertible in $\mathcal{C}_0(X).$ From Lemma \ref{idempotent},  there exists two idempotents $P , Q  \in L(X)$ such that $\pi(P) \in comm(\pi(T))$, $\pi(Q) \in comm(\pi(T)),$  $\pi(T)$ is left $\pi(P)-$invertible in $\mathcal{C}_0(X),$ $\pi(T)\pi(P) $ is nilpotent in $\mathcal{C}_0(X),$   $\pi(T)$ is right $\pi(Q)-$invertible in $\mathcal{C}_0(X)$ and $\pi(T)\pi(Q) $ is nilpotent in $\mathcal{C}_0(X).$ As $\pi(T)$ is left $\pi(P)-$invertible in $\mathcal{C}_0(X),$ it has a left inverse  $\pi(S), S\in L(X)$ such that $\pi(S)\pi(P)= \pi(P) \pi(S).$ Then  $\pi(I-P)\pi(T)\pi(I-P)$ is left invertible in the algebra $\pi(I-P)\mathcal{C}_0(X)\pi(I-P)$ whose identity element is $ \pi(I-P), $  having as left inverse $\pi(I-P) \pi(S)\pi(I-P)$ which commutes with $\pi(P).$
So $$\pi(I-P) \pi(S)\pi(I-P) \pi(I-P)\pi(T)\pi(I-P)= \pi(I-P).$$
Similarly, as $\pi(T)$ is right $\pi(Q)-$invertible in $\mathcal{C}_0(X),$ it has a right inverse  $\pi(R), R\in L(X)$ such that $\pi(R)\pi(Q)= \pi(Q) \pi(R).$ Then  $\pi(I-Q)\pi(T)\pi(I-Q)$ is right invertible in the algebra $\pi(I-Q)\mathcal{C}_0(X)\pi(I-Q)$ whose identity element is $ \pi(I-Q), $  having as right inverse $\pi(I-Q) \pi(R)\pi(I-Q)$ which commutes with $\pi(Q).$
So $$  \pi(I-Q)\pi(T)\pi(I-Q) \pi(I-Q) \pi(R)\pi(I-Q)= \pi(I-Q).$$

For $ n \in \mathbb{N},$  we have $ \pi(T)^n= (\pi(P)\pi(T)\pi(P) +   \pi(I-P)\pi(T)\pi(I-P))^n= (\pi(Q)\pi(T)\pi(Q) +   \pi(I-Q)\pi(T)\pi(I-Q))^n.$ Hence for n large enough, $ \pi(T)^n=  \pi(I-P)\pi(T)^n\pi(I-P)=   \pi(I-Q)\pi(T)^n\pi(I-Q),$ because $\pi(P)\pi(T)\pi(P)$ and $\pi(Q)\pi(T)\pi(Q)$ are nilpotent. 

Then $ 0= \pi(P) \pi(T)^n= \pi(P) \pi(I-Q) \pi(T)^n  \pi(I-Q).$ Thus $  \pi(P) \pi(I-Q)= \pi(P) \pi(I-Q) \pi(T)^n \pi(I-Q) \pi(I-Q)\pi(R)^n\pi(I-Q)=0. $ Therefore $\pi(P) =\pi(PQ).$

Similarly we have $0= \pi(T)^n \pi(Q)= \pi(I-P)\pi(T)^n\pi(I-P)\pi(Q).$ Thus $\pi(I-P)\pi(Q)=  \pi(I-P) \pi(S)^n \pi(I-P) \pi(I-P)\pi(T)^n\pi(I-P)\pi(Q)= 0.$ Therefore $\pi(Q) =\pi(PQ).$
 and $\pi(T)$ is  Drazin invertible in $\mathcal{C}_0(X).$  \ep

\bet \label{left-Drazin} Let $T \in L(X).$ Then the following properties
are equivalent:

\begin{enumerate}

\item $T$ is a left semi-B-Fredholm operator.
\item  $T$ is a nilpotent\textbf{\Large-}left semi-Fredholm operator

\end{enumerate}
\eet

\bp  $1)\Rightarrow 2)$  Assume that $T$ is a left semi-B-Fredholm operator, so  $\pi(T)$ is  left  Drazin invertible in  $\mathcal{C}_0(X).$ Then there exist an
idempotent $p \in \mathcal{C}_0(X)$ so that:  

  $\bullet$    $p \pi(T) = \pi(T) p.$

   $\bullet$ $p \pi(T)$ is nilpotent in $\mathcal{C}_0(X).$

   $\bullet$  There exists  $S \in L(X)$  such that $p \pi(S) = \pi(S) p$ and  $\pi(S) ( \pi(T) + p)) = \pi(I).$


From Lemma \ref{idempotent}, there exists an idempotent $ P \in L(X)$ such that $ \pi(P)=p.$

 Since $\pi(T)$ and $\pi(P)$ commutes, we have  $\pi(PTP)=\pi(TP)$ and  it follows that $PTP$
 is a power finite rank operator. Moreover from Remark \ref{not1}, we have

 \begin{equation}\label{Decomposition}
 T=T_1\oplus T_2 +F,
 \end{equation}

   where $F$ is a finite rank operator.

 Let us show that $T_2$ is a left semi-Fredholm operator. We have  $\pi(S) ( \pi(T) + \pi(P))) = \pi(I)$  and  $\pi(T) \pi(I-P)=\pi(I-P)\pi(T),$ so  $ (I-P)S(I-P) (I-P) (T+P)(I-P)= I-P  + (I-P)F'(I-P)$  and $ [(I-P)S(I-P)]_{\mid{X_2}}  T_2= I_2  + [(I-P)F'(I-P)]_{\mid X_2},$
where   $F'$ is a  finite rank operator.
Hence $T_2$ is a left semi-Fredholm operator, because $\pi(T_2)$ is left invertible  in the algebra $\mathcal{C}_0(X_2).$

If $n$ is large enough, then $ (T_1)^n$ is a finite rank operator. So the operator $(T_1)_{[n]}: R(T_1^n) \rightarrow R(T_1^n)$ is a Fredholm operator. Thus $T_1$ is a B-Fredholm operator and from \cite[Proposition2.6]{P7}, $ T_1$ is a quasi-Fredholm  operator. Then obviously, $T_1 \oplus T_2$ is a quasi-Fredholm operator. 

Let $S= T_1 \oplus T_2$ and let $d=dis(S).$ Then  we have $ N(S)\cap R(S^d)= [N(T_1)\cap R(T_1^n)] \oplus [N(T_2)\cap R(T_2^d)]$ and   $ R(S)+ N(S^d)= [R(T_1)+ N(T_1^d)] \oplus [R(T_2)+ N(T_2^d)].$ 

As $S$ is a quasi-Fredholm operator, then from  \cite[Th\'eor\`eme 3.2.2]{LAB}   there exists two closed
subspaces  $M,  N$  of $ X $ such that   $ X = M\oplus N $  and

\begin{enumerate}
\item $S(N)\subset N $  and  $ S_{\mid N}$ is a  nilpotent operator of degree $d$,
\item $ S(M) \subset M,     
 N(S_{\mid M })  \subset \underset{m} \cap R( (S_{\mid M })^m)  $  and $S_{\mid M }$ is a regular operator.
\end{enumerate}

$\bullet$ It is easily seen that  $ N(S_{\mid M })= N(S)\cap R(S^d)= (N(T_1)\cap R(T_1^d)) \oplus (N(T_2)\cap R(T_2^d)).$  
As $T_1$ is a B-Fredholm operator and $ T_2$ is a left semi-Fredholm operator, then $ N(S_{\mid M })$ is of finite dimension.

$\bullet$ As $S_{\mid M }$ is a regular operator, then $R(S_{\mid M })$ is complemented in $M.$ 

$\bullet$ Thus $S$ is a nilpotent \textbf{-} left semi-Fredholm operator. 

Now as $ T= S+ F,$ with $F$ of finite rank, then from Theorem \ref{perturbation},   $ T$ is a quasi-Fredholm operator. So from  \cite[Th\'eor\'eme 3.2.2]{LAB}   there exists two closed
subspaces  $M',  N'$  of $ X $ and an integer $d'= dis(T)$ such that   $ X = M'\oplus N' $  and

\begin{enumerate}
\item $T(N')\subset N' $  and  $ T_{\mid N'}$ is a  nilpotent operator of degree $d'$,
\item $ T(M) \subset M,  
 N(T_{\mid M' })  \subset \underset{m} \cap R( (T_{\mid M' })^m)  $  and $T_{\mid M' }$ is a regular operator 
\end{enumerate}

$\bullet$ Similarly to the case of $S$ and using \cite{KMMP},  we  have $ N(T_{\mid M' })= N(T)\cap R(T^{d'})\overset{e}{=} N(S)\cap R(S^d).$ As  $N(S)\cap R(S^d)$  is of finite dimension,  then from Theorem \ref{perturbation}, $ N(T_{\mid M' })$ is of finite dimension.

$\bullet$ As $T_{\mid M' }$ is a regular operator, then $R(T_{\mid M' })$ is complemented in $M'.$ 

$\bullet$ Thus $T$ is a nilpotent \textbf{-} left semi-Fredholm operator.

$ 2)\Rightarrow 1)$ Conversely assume that there exists an idempotent $P \in L(X)$ such that $PT= TP,$   $T_1 $ is a  nilpotent  operator  and  $T_2$ is a left semi-Fredholm one. We have $ T= T_1 \oplus T_2.$ 

So  $ T_1\oplus T_2 +P = ( T_1  \oplus 0) +  P  + ( 0 \oplus T_2) = P[ (T_1 \oplus 0)+ I)] P  + (I-P) ( I_1 \oplus  T_2 ) (I-P).$  As $T_1$ is a nilpotent operator, then $T_1 \oplus 0$ is also a nilpotent  operator and $\pi( (T_1 \oplus 0)+ I))=  \pi( (T_1 \oplus 0)) + \pi(I)$ is invertible in $\mathcal{C}_0(X).$ Let $ \pi(S_1) $ be its inverse, where $ S_1 \in L(X).$

As $T_2$ is a left  semi-Fredholm operator in $L(X_2),$  then from \cite[Corollary 2.4]{HB2}, there exists $ S_2 \in  L(X_2)$  such that $ S_2 T_2- I_2$ is a finite rank operator. Moreover $I_1 \oplus S_2$ commutes with $P$   because
$(I_1  \oplus S_2)P=P(I_1 \oplus S_2)=I_1 \oplus 0.$
We observe that $\pi( I_1\oplus  T_2) $ is left invertible in  $\mathcal{C}_0(X)$ having  $ \pi( I_1 \oplus S_2)$ as a left inverse. Then:
$$ \pi((PS_1P + (I-P)(I_1\oplus S_2) (I-P))) \pi ( T+P)$$ $$=\pi((PS_1P + (I-P)(I_1\oplus S_2) (I-P))) \pi (T_1\oplus T_2 +P)$$
$$=\pi([PS_1P + (I-P)(I_1\oplus S_2) (I-P)]) \pi ( P[ (T_1 \oplus 0)+ I)] P  + (I-P) ( I_1 \oplus  T_2 ) (I-P))=\pi(I).$$
It is easily seen that  $\pi((PS_1P + (I-P)( I_1 \oplus S_2) (I-P)))$ commutes with $\pi(P).$

\noindent Moreover we have $\pi(P)\pi(T)= \pi( T_1 \oplus 0)$   is  nilpotent in $\mathcal{C}_0(X)$  because $ T_1 \oplus 0$ is a nilpotent operator.  Thus $ T$ is a left  Drazin invertible in $\mathcal{C}_0(X).$  \ep


 As the  proof of the next result is very similar to the proof of  Theorem \ref{left-Drazin}, we  include it  under a  lightened version.

\bet \label{right-pseudo-0} Let $T \in L(X).$ Then the following properties
are equivalent:

\begin{enumerate}

\item $T$ is a right semi-B-Fredholm operator.
\item $T$  is  a nilpotent \textbf{-} right semi-Fredholm operator.

\end{enumerate}

\eet

\bp  $1)\Rightarrow 2)$  Assume that $T$ is a right semi-B-Fredholm operator. As in the proof of Theorem \ref{left-Drazin} there exist an
idempotent $P \in L(X),$ such  that   $ X= P(X) \oplus (I-P)(X)= X_1 \oplus X_2$ and relatively to this decomposition, we have 
 \begin{equation}\label{Decomposition}
 T=T_1\oplus T_2 +F,
 \end{equation}

where  $T_1$ is a  B-Fredholm operator, $T_2$ is a right semi-Fredholm operator and  $F$ is a finite rank operator.

Let $S= T_1 \oplus T_2$ and let $d=dis(S),$ then  we have $ N(S)\cap R(S^d)= N(T_1)\cap R(T_1^n) \oplus N(T_2)\cap R(T_2^d)$ and   $ R(S)+ N(S^d)= (R(T_1)+ N(T_1^d)) \oplus (R(T_2)+ N(T_2^d)).$ 

As $S$ is a quasi-Fredholm operator, then from  \cite[Th\'eor\`eme 3.2.2]{LAB}   there exists two closed
subspaces  $M,  N$  of $ X $  such that   $ X = M\oplus N $  and:

\begin{enumerate}
\item $S(N)\subset N $  and  $ S_{\mid N}$ is a  nilpotent operator of degree $d$,
\item $ S(M) \subset M, 
N(S_{\mid M })  \subset \underset{m} \cap R( (S_{\mid M })^m)  $  and $S_{\mid M }$ is a regular operator
\end{enumerate}

$\bullet$ As $S_{\mid M }$ is a regular operator, then  $N(S_{\mid M })$ is complemented in $M.$

$\bullet$ Similarly  $R(S_{\mid M }) \oplus N=  R(S)+ N(S^d)= (R(T_1)+ N(T_1^d)) \oplus (R(T_2)+ N(T_2^d))$ is of finite codimension in $X,$  because $T_1$ is a B-Fredholm operator and $T_2$ is a left semi-Fredholm operator. Indeed $ (R(T_1)+ N(T_1^d))$ is of finite codimension in $ X_1=P(X)$ and $(R(T_2)+ N(T_2^d))$ is of finite codimension in $X_2=(I-P)(X_2).$
So there exists a finite dimensional  subspace $E$ of $X$ such that  $ X= X_1 \oplus X_2= R(S_{\mid M }) \oplus N \oplus E= R(S_{\mid M }) \oplus (N \oplus E) .$  Then  $ M= R(S_{\mid M }) \oplus [(N \oplus E)\cap M]= R(S_{\mid M }) + P'(E), $ where $ P': X \rightarrow X$ is the projection on $M$ in parallel to $N.$   Then  $R(S_{\mid M })$ is of finite codimension  in $M,$  because $E$ 
is of finite dimension.

$\bullet$ Thus $S$ is a nilpotent\textbf{\Large-}left semi-Fredholm operator. 

Now as $ T= S+ F,$ with $F$ of finite rank, then from Theorem \ref{perturbation}, $ T$ is a quasi-Fredholm operator. So from  \cite[Th\'eor\'eme 3.2.2]{LAB}   there exists two closed
subspaces  $M',  N'$  of $ X $ and an integer  such that   $ X = M'\oplus N' $  and:

\begin{enumerate}
\item $T(N')\subset N' $  and  $ T_{\mid N'}$ is a  nilpotent operator of degree $d'$,
\item $ T(M') \subset M',  
 N(T_{\mid M' })  \subset \underset{m} \cap R( (T_{\mid M' })^m)  $ and $T_{\mid M' }$ is a regular operator
\end{enumerate}

$\bullet$ As $T_{\mid M' }$ is a regular operator, then  $N(T_{\mid M' })$ is complemented in $M'.$

$\bullet$ Moreover  $R(T_{\mid M' }) \oplus N'=  R(T) + N(T^{d'}) \overset{e}{=} R(S) + N(S^d).$ As $ R(S) + N(S^d)$ is of finite codimension in $X,$   then  $R(T_{\mid M' }) \oplus N'$ is also of finite codimension in $X.$ Then with the same method as in the case of $S,$  we can see that  $R(T_{\mid M' })$ is of finite codimension  in $M'.$

$\bullet$ Thus $T$ is a nilpotent \textbf{-} right semi-Fredholm operator.

$ 2)\Rightarrow 1)$ We follow  the same method as in $ 2)\Rightarrow 1)$ in  Theorem \ref{left-Drazin}  \ep

\noindent  Alaa Hamdan\\
\noindent   Dubai \\
\noindent   United Arab Emirates \\
\noindent   \href{mailto:aa.hamdan@outlook.com }{aa.hamdan@outlook.com}\\

\noindent  Mohammed Berkani\\
\noindent Honorary member of LIAB \\
\noindent Department of Mathematics\\
\noindent Science Faculty of Oujda \\
\noindent Mohammed I University\\
\noindent Morocco\\
\noindent \href{mailto:berkanimoha@gmail.com}{berkanimoha@gmail.com}

\end{document}